\documentclass[12pt,twoside]{article}
\usepackage{amsxtra,amssymb,amsthm,amsmath,latexsym}

\theoremstyle{plain}

\def\oH{\buildrel\circ\over H}
\def\oH1{\buildrel\circ\over H\kern-.02in{}^1}
\def\oH1{\buildrel\circ\over H\kern-.02in{}^1}

\def\d{\delta}

\def\a{\alpha}

\begin{document}

\title{On unbounded operators and applications
\thanks{Math subject classification: 47A05, 45A50,35R30}
\thanks{key words:  unbounded linear operators, 
ill-posed problems, regularization, 
discrepancy principle }
}

\author{                          
A.G. Ramm\\       
 Mathematics Department, Kansas State University, \\
 Manhattan, KS 66506-2602, USA\\
ramm@math.ksu.edu\\}

\date{}

\maketitle \thispagestyle{empty}

{\bf Abstract} 

Assume that $Au=f,\quad (1)$
is a solvable linear equation in a Hilbert space $H$, $A$ is a linear, 
closed, densely defined, unbounded operator in $H$, which is not boundedly 
invertible, so problem (1) is ill-posed. 
It is proved that the closure of the operator $(A^*A+\a I)^{-1}A^*$, with 
the domain $D(A^*)$, where $\a>0$ is a constant, is a linear
bounded everywhere defined operator with norm $\leq 1$.
This result is applied to the variational problem $F(u):=
||Au-f||^2+\a ||u||^2=min$, where $f$ is an arbitrary element of $H$,
not necessarily belonging to the range of $A$.
Variational regularization of problem (1) is constructed, and a 
discrepancy principle is proved.

1. {\bf Introduction}

 The main results of this paper are formulated as 
Theorems 1 and 2
and proved in Sections 1 and 3, respectively.
In Section 1 we formulate Theorem 1 which deals with a linear, unbounded, 
closed, densely defined operator $A$. In Section 2 this operator is
assumed not boundedly invertible and the problems arising in
the study of variational regularization of the solution to the equation
$$Au=f,
 \eqno{(1)}
$$ 
are studied, where $A:H\to H$
is a linear, unbounded, closed, densely defined, not boundedly invertible 
operator on a 
Hilbert space $H$ with domain $D(A)$ and range $R(A)$. Since $A$ is
densely defined and closed, its  adjoint 
$A^*$ is a closed, densely defined 
linear operator. The operators $T=A^*A$ and $Q=AA^*$ are nonnegative, 
selfadjoint,
densely defined in $H$ operators (see [1]),  the operator $T_\a:=T+\a I$,
($I$ is the identity operator and $\a>0$ is a constant) is boundedly 
invertible, i.e.,
its inverse is a bounded linear operator, defined on all of $H$,
with norm $\leq \frac 1 \a$.
It is easy to check  that the operator $A^*Q^{-1}_\a$ is bounded, 
defined on all of $H$, and $||A^*Q^{-1}_\a||\leq \frac 1{2\sqrt{\a}}$. 
We assume in Section 2 that the operator $A$ is not boundedly 
invertible, in which case problem (1) is ill-posed.

We are interested in the operator $S:=S_a:=T^{-1}_\a A^*$ defined on
a dense set $D(A^*)$.   The reasons for our interest will be explained soon. 
The product of an unbounded closed operator ($A^*$ in our case) and a 
bounded operator ($T^{-1}_\a$ in our case) is not necessarily 
closed, in general, 
as a simple example shows: if $A=A^*\geq 0$ is unbounded 
selfadjoint (and consequently closed) operator and $B=(I+A)^{-1}$ is a 
bounded operator, then
the operator $BA$ with the domain  $D(A)$ is not closed. Its closure is
a bounded operator defined on all of $H$. This closure is uniquely defined 
by continuity.

{\bf Lemma 1.} {\it If $A$ is a linear, closed, densely defined, unbounded 
operator in $H$, and $B$ is a bounded linear operator such that
$R(B^*)\subset D(A)$, then the operator $BA^*$ with the domain $D(A^*)$
is closable.}

{\bf Proof.} To prove the closability of $BA^*$ one has to 
prove that if
$u_n\to 0$ and $BA^*u_n\to w$, then $w=0$. Let $h$ be arbitrary. Then
$B^*h$ belongs to $D(A)$. Therefore
$$(w,h)=\lim (BA^*u_n, h)=\lim (u_n, AB^*h)=0.$$
Thus, $w=0$. \hfill $\Box$

In our case $B=T^{-1}_\a=B^*$ and $R(T^{-1}_\a)\subset D(A)$. 
By Lemma 1, 
the operator $T^{-1}_\a A^*$ is closable. The operator 
$AT^{-1}_\a$ is bounded, defined on all of $H$,  with norm 
$\leq \frac 1 {2\sqrt{\a}}$. Indeed, by the polar decomposition
one has $A=UT^{1/2}$, where $U$ is an isometry, so $||U||\leq 1$.
Thus, $||AT^{-1}_\a||\leq ||T^{1/2}T^{-1}_\a||=\sup_{s\geq 0}\frac 
{s^{1/2}}{s+\a}= \frac 1 {2\sqrt{\a}}$. 

{\bf Lemma 2.} {\it The operator $S:=T^{-1}_\a A^*$ with domain $D( A^*)$
has the closure $\overline{S}=S^{\ast \ast}$, which is a bounded operator 
defined on all of $H$, with norm $\leq \frac 1 {2\sqrt{\a}}$, and
$S^*=AT^{-1}_\a$ is a bounded operator defined on all of $H$.} 
  
{\bf Proof.} The operator $S$ is densely defined. By Lemma 1 it is 
closable, so the operator $S^*$ is densely defined.
Let us prove that $S^*=AT^{-1}_\a$. Let $h\in H$ be arbitrary. We have
$$(T^{-1}_\a A^*u,h)=(A^*u, T^{-1}_\a h)=(u, AT^{-1}_\a h).$$
This implies that $D(S^*)=H$ and $(T^{-1}_\a A^*)^*=AT^{-1}_\a$. We have 
used the relation 
$R(T^{-1}_\a) \subset D(A)$. Let us check this relation. Let 
$g\in R(T^{-1}_\a)$, then $g=T^{-1}_\a h$ and $h=T_\a g$. Thus,
$g\in D(T)\subset D(A)$, as claimed.
Lemma 2 is proved. \hfill $\Box$

From Lemmas 1 and 2 we obtain the following result.

{\bf Theorem 1.} {\it Let $A$ be a linear, closed, densely defined, 
unbounded operator in a Hilbert space.
Then the operator $S=T^{-1}_\a A^*$ with domain $D(A^*)$
admits a unique closed extension $\overline {S}$ defined on all of $H$,
with the norm $\leq \frac 1{2\sqrt{\a}}$.} 
\vskip .07in
{\it Why should one be interested in the above theorem?} 

The answer is: because of its crucial role in the study 
of equation (1) and of variational 
regularization for equation (1). The corresponding theory is developed in 
Section 2 and proofs are given in Section 3.

2. {\bf  Variational regularization}

{\bf Assumption A}: {\it We assume throughout that $A$ is linear, 
unbounded, densely defined 
operator in $H$, and that $A$ is not boundedly invertible, so problem (1) 
is ill-posed. We assume that equation (1) is solvable, possibly 
nonuniquely, that $f\neq 0$, and denote by $y$ its unique 
minimal-norm solution,
$y\bot N$, where $N=N(A)$.} 

This assumption is not repeated below, but is a standing one throughout 
the rest of this paper.

Assume that $||f_\d-f||\leq \d$, where $f_\d$ 
is the "noisy" data, which are known for some given small $\d>0$, while 
the 
exact data $f$ are unknown. The problem is to construct a stable 
approximation $u_\d$ to $y$, given the data $\{A, \d, f_\d\}$.
Stable approximation means that $\lim_{\d \to 0}||u_\d-y||=0$.
Variational regularization is one of the methods for constructing such an
approximation.

If $A$ is bounded, this method consists of 
solving the minimization problem
$$
F(u):=F_{\a,\d}=||Au-f_\d||^2 +\a||u||^2=\min, 
 \eqno{(2)}
$$
and choosing the regularization parameter $\a=\a(\d)$ so that
$\lim_{\d\to 0}u_{\d}=y$, where $u_{\d}:=u_{\a(\d),\d}$.
It is well known and easy to prove that if $A$ is bounded, 
then
problem (2) has a unique solution,
$u_{\a,\d}= T^{-1}_\a A^*f_\d$, which is a unique global minimizer
of the quadratic functional (2), and this minimizer solves the 
equation $T_\a u_{\a,\d}=A^* f_\d$. The last equation does not make sense,
in general, if $A$ is unbounded, because $f_\d$ may not belong to 
$D(A^*)$. This is the difficulty arising in the case of unbounded $A$.
In this case it is not a priori clear if the global minimizer of
functional (2) exists. We prove that this minimizer exists for any $f_\d 
\in H$, that it is unique, and that there is a function $\a=\a(\d)>0$,
$\lim_{\d\to 0}\a(\d)=0$, such that $\lim_{\d\to 0}u_{\a(\d),\d}=y$, so
the element $u_\d:=u_{\a(\d),\d}$ is a stable approximation 
of the 
unique minimal-norm solution to (1). Theorem 1 allows one to define the 
element $T^{-1}_\a A^*f_\d$ for any $f_\d$, and not only for those $f_\d$
which belong to $D(A^*)$.
We also prove for unbounded $A$ a discrepancy principle 
in the following form.
Let $u_{\d, \a}$ solve (2). Consider the  equation
for finding $\a=\a(\d)$:
$$ ||Au_{\a,\d}-f_\d||=C\d,\quad C>1, \quad ||f_\d||>C\d,
 \eqno{(3)}
$$
where $C$ is a constant. Equation (3) is the discrepancy 
principle. We prove that equation (3) determines 
 $\a(\d)$  uniquely, $\a(\d)\to 0$ as $\d\to 0$,
and $u_\d:=u_{\d, \a(\d)}\to y$ as $\d\to 0$. This justifies 
the discrepancy principle for choosing the regularization 
parameter (see [2] for various forms of the discrepancy principle).

Let us formulate the results.

{\bf Theorem 2.} {\it For any $f\in H$ the functional  
$F(u)=||Au-f||^2 +\a ||u||^2$ has a unique global minimizer
$u_\a= A^*Q_\a^{-1}f$, where $Q=AA^*$, $Q_\a:=Q+\a I$,  $\a>0$ is a 
constant, and 
$$A^*Q_\a^{-1}f=T^{-1}_\a A^*f,
\eqno{(4)}
$$
where $T^{-1}_\a A^*$ is the closure of the operator $T^{-1}_\a A^*$
defined on $D(A^*)$.  If $f\in R(A)$, then
$$\lim_{\a \to 0}||u_\a-f||=0,
\eqno{(5)}
$$
where $u_\a$ is the unique global minimizer of $F(u)$ and $y$ is the 
minimal-norm solution to (1).
If $||f_\d -f||\leq \d$ and $u_{\a,\d}=A^*Q_\a^{-1}f_\d$,
then there exists an $\a(\d)>0$ such that
$$\lim_{\d \to 0}||u_\d-f||=0,\quad \lim_{\d \to 0}\a(\d)=0, \quad u_\d:=
u_{\a(\d),\d}.
\eqno{(6)}
$$
Equation (3) is uniquely solvable for $\a$, and for its solution $\a(\d)$
equation (6) holds.
} 

In Section 3 proof of Theorem 2 is given.

3. {\bf Proof of Theorem 2}

3.1. For any $h\in D(A)$ let $u_\a:=A^* Q_\a^{-1}$. One has
$$ F(u_\a+h)=F(u_\a) +||Ah||^2 +\a ||h||^2 +
2 \Re [(Au_\a-f,Ah)+\a(u_\a,h)],
\eqno{(7)}
$$
and 
$$(Au_\a -f,Ah)+\a(u_\a,h)=(Q 
Q_\a^{-1}f-f,Ah)+\a(u_\a,h)=-\a(Q_\a^{-1}f,Ah)+\a(u_\a,h)=
$$
$$-\a[(A^* Q_\a^{-1}f+u_\a,h)=0.
\eqno{(8)}
$$
From (7) and (8) it follows that $u_\a$ is the 
unique global minimizer of $F(u)$. 


3.2. Let us prove (4). If (4) holds on a dense in $H$ linear 
subset $D(A^*)$, then it holds on all of  $H$ by continuity 
because
$A^* Q_\a^{-1}$ is a bounded linear operator, defined on all 
of $H$, with norm $\leq \frac 1 {2\sqrt{\a}}$, so that the 
closure of the operator $T_\a^{-1}A^*$ defined on $D(A^*)$,
is a bounded 
operator, defined on all of $H$, with norm $\leq \frac 1 
{2\sqrt{\a}}$.   Indeed, let 
$f\in D(A^*)$, $g:= Q_\a^{-1}f$, so $ Q_\a g=f$
and $g\in D(A^* A A^*)$. Therefore
equation (4) is equivalent to
$A^*Q_\a g=T_\a A^* g$, or $A^*AA^*g+\a A^*g=A^*AA^*g+\a 
A^*g$, which is an identity. If $f\in D(A^*)$ (so that
$g\in D(A^*AA^*)$), then the above formulas are justified
and one can go back from the identity
$A^*AA^*g+\a A^*g=A^*AA^*g+\a A^*g$, valid for any $g\in  
D(A^*AA^*)$, define $f=Q_\a g$, (this $f$ belongs to 
$ D(A^*)$ because $Qg\in  D(A^*)$), and get (4).

Note that if $A$ were bounded, then one would have 
the identity 
$$
A^*\phi (Q)=\phi(T)A^*, \quad T=A^*A, \quad Q=AA^*,
\eqno{(9)}
$$
valid for any continuous function $\phi$. Indeed, if $\phi$
is a polynomial, then (9) is obvious (for example,
if $\phi(s)=s$, then (9) becomes $A^*(AA^*)=(A^*A)A^*$). 
If $\phi$ is a continuous function on the interval $[0, 
||A||^2]$, then it is a limit
(in the $\sup$-norm on this bounded interval) 
of a sequence of polynomials (Weierstrass' theorem), so (9) 
holds. In our problem $A$ is unbounded, so are $Q$ and 
$T$, and $\phi (s) =\frac 1 {s+\a}$ (with 
$\a=const>0$) is a continuous function on an infinite 
interval $[0,\infty)$. Linear unbounded operators do not 
form an algebra, in general, because of the difficulties 
with domain of definition of the product of two unbounded 
operators (the product may have the trivial domain 
$\{0\}$). That is why formula (4), which is a particular 
case of (9) for bounded operators, has to be proved
independently of this formula.

3.3 Let us prove (5). If $f\in R(A)$, then $f=Ay$, where 
$y\bot N$ is the minimal-norm solution to (1). We have
$u_\a-y=T^{-1}_\a Ty-y=-\a T^{-1}_\a y$ and
$$
\lim_{\a\to 0}||\a T^{-1}_\a y||^2=\lim_{\a\to 
0}\int_0^\infty \frac 
{\a^2}{(\a+s)^2}d(E_sy,y)=||P_Ny||^2=0,
$$
where $E_s$ is the resolution of the identity of the 
selfadjoint operator $T$ and $P_N$ is the orthogonal 
projector onto $N=N(A)$, so $||P_Ny||=0$ because $y\bot N$.

3.4. Let us prove (6). We have 
$$ ||u_\d-y||\leq ||u_\d-u_\a||+||u_\a-y||:=I_1 +I_2.$$
We have already proved that $\lim_{\d\to 0}I_2=0$,
because $\lim_{\d\to 0}\a(\d)=0$. Let us estimate $I_1$:
$$ ||u_\d-u_\a||=||A^*Q_{\a(\d)}^{-1}(f_\d-f)||\leq \frac \d 
{2\sqrt{\a}}$$  
Thus, if $\lim_{\d\to 0}\a(\d)=0$ and $\lim_{\d\to 0}\frac 
\d{2\sqrt{\a}}=0$, then (6) holds.

3.5. Finally, let us prove 

{\it The discrepancy principle:

Equation (3) is uniquely solvable for $\a$ and its solution 
$\a(\d)$  satisfies (6).}

The proof follows the one in [2], p.22.
One has 
$$g(\a,\d):= 
||Au_{\a,\d}-f_\d||^2=||QQ_\a^{-1}f_\d-f_\d||^2=
\a^2\int_0^\infty
\frac {d(\mathcal{E}_s f_\d,  f_\d)}{(s+\a)^2}=C^2\d^2,
\eqno{(10)}
$$
where $\mathcal{E}_s$ is the resolution of the identity
of the selfadjoint operator $Q$. The function 
$g(\a):=g(\a,\d)$ for a fixed $\d>0$ is continuous, 
strictly increasing on $[0,\infty)$ and 
$g(\infty)>C^2\d^2$ while $g(0)\leq \d^2$, as we will prove
below. Thus,  there exists a unique $\a(\d)>0$, such that
$g(\a(\d),\d)=C^2\d^2$, and $\lim_{\d \to 0}\a(\d)=0$ 
because $g(\a,\d)>0$ for $\a\neq 0$ and any $\d\in 
[0,\d_0)$, provided that $||f||\neq 0$, which we assume. 
Here $\d_0>0$ is a 
sufficiently small number.

Let us prove the two inequalities: $g(\infty)>C^2\d^2$
and $g(0)\leq \d^2$. We have
$$g(\infty)=\int_0^\infty d(\mathcal{E}_s f_\d,  
f_\d)=||f_\d||^2>C^2\d^2,$$ 
because of the assumption 
$||f_\d||>C\d$. Also
$$g(0)=||P_{N(Q)}f_\d||^2\leq \d^2.$$
Indeed, $f_\d=f+(f_\d-f)$. The element $f\in R(A)$, so
$f\bot N(A^*)=N(Q)$. Therefore 
$P_{N(Q)}f_\d=P_{N(Q)}(f_\d-f)$. Consequently,
$$||P_{N(Q)}f_\d||\leq ||P_{N(Q)}(f_\d-f)||\leq 
||f_\d-f||\leq \d.$$ 

Let us prove the limiting relation $\lim_{\d \to 
0}||u_\d-f||=0$. We have
$$
F_{\a(\d),\d}(u_\d)=||Au_\d -f_\d||^2 +\a(\d)||u_\d||^2\leq
F_{\a(\d),\d}(y)\leq \d^2+\a(\d)||y||^2.
\eqno{(11)}
$$
Since $||Au_\d -f_\d||^2=C^2\d^2>\d^2$ we conclude from (11) 
that $||u_\d||\leq ||y||$ for all $\d\in [0, \d_0)$.
Thus we may assume that $u_\d\rightharpoonup z $ as $\d\to 
0$, where $\rightharpoonup$ denotes weak convergence in $H$.
Since $\lim_{\d\to 0}f_\d=f$, 
we conclude from  $||Au_\d -f_\d||=C\d$ that $\lim_{\d\to 
0}||Au_\d-f||=0$.
This implies  $Az=f$. Indeed, for any $h\in D(A^*)$ 
one has 
$$(f,h)=\lim_{\d\to 0}(Au_\d,h)=(z,A^*h).$$
 Therefore
$Az=f$. Since $||u_\d||\leq ||y||$, we have 
$\overline{\lim}_{\d\to 0}||u_\d||\leq ||y||$.
From $u_\d \rightharpoonup z$ we obtain $||z||\leq 
\underline{\lim}_{\d\to 0}||u_\d||\leq ||y||$.
Thus, $||z||\leq ||y||$. Since the minimal-norm solution to 
(1) is unique, it follows that $z=y$. Thus, 
$u_\d \rightharpoonup y$ and $\overline{\lim}_{\d\to 
0}||u_\d||\leq ||y||\leq \underline{\lim}_{\d\to 
0}||u_\d||$. This implies $\lim_{\d\to 0}||u_\d||=||y||$.
Consequently, $\lim_{\d\to 0}||u_\d-y||=0$, because
$||u_\d-y||^2=||u_\d||^2+||y||^2-2\Re (u_\d,y)\to 0$ as 
$\d\to 0$.

Theorem 2 is proved. \hfill $\Box$

\end{document}